\def\mytitle{Dual-Quaternion Fourier Transform}

\def\mysubtitle{}

\def\mykeywords{Dual-Quaternions; Signals; Fourier; Dual-Quaternion Fourier Transform; Analytics; Quaternions; Dual-Numbers}

\def\myauthor{Benjamin Kenwright}

\def\myabstract{
Fourier transform (FT) plays a crucial role in a broad range of applications, from enhancement, restoration and analysis through to security, compression and manipulation. The Fourier transform (FT) is a process that converts a function into a form that describes the frequencies. This process has been extended to many domains and numerical representations (including quaternions). However, in this article, we present a new approach using dual-quaternions. As dual-quaternions offer an efficient and compact symbolic form with unique mathematical properties. While dual-quaternions have established themselves in many fields of science and computing as an efficient mathematical model for providing an unambiguous, un-cumbersome, computationally effective means of representing multi-component data, not much research has been done to combine them with Fourier processes. Dual-quaternions are simply the unification of dual-number theory with hypercomplex numbers; a mathematical concept that allows multi-variable data sets to be transformed, combined, manipulated and interpolated in a unified non-linear manner. We define a Dual-Quaternion Fourier transform (DQFT) for dual-quaternion valued data over dual-quaternion domains. This opens the door to new types of analysis, manipulation and filtering techniques.  We also present the Inverse Dual-Quaternion Fourier Transform (IDQFT). The DQFT unlocks the potential provided by hypercomplex algebra in higher dimensions useful for solving dual-quaternion partial differential equations or functional equations (e.g., for multicomponent data analysis). 
}

\documentclass[9pt,journal,cspaper,compsoc]{IEEEtran}

\usepackage{multicol}

\usepackage[backref=page]{hyperref}

\hypersetup{
	colorlinks = true, 
	linkcolor = blue,
	anchorcolor = red,
	citecolor = blue, 
	filecolor = red, 
}

\hypersetup{pdfinfo={
		Author   		= {\myauthor},
		Title    		= {\mytitle},
		Subject  		= {\mytitle},
		Keywords 		= {\mykeywords},
		CreationDate 	= {D:2014020220195600},
		Producer		= {\myauthor},
		Creator			= {\myauthor},
}}

\usepackage{bookmark}

\usepackage{multicol}

\usepackage{url}

\usepackage{atbegshi,picture}

\usepackage{tikz}
\usepackage{calc}% http://ctan.org/pkg/calc
\usepackage{setspace}

\usetikzlibrary{mindmap,shadows}

\usepackage{pgfplots}

\usepackage{pgf-pie}

\usepackage{smartdiagram}

\usepackage{bchart}

\usepackage{amsmath} % /checkmark
\usepackage{amssymb}

\renewcommand{\theparagraph}{}
\usepackage{titlesec}
\titleformat{\paragraph}[runin]{\normalfont\normalsize\bfseries}{\theparagraph}{0em}{}
\titlespacing*{\paragraph}{0pt}{1.0ex plus 2ex minus .1ex}{1em}

\titlespacing*{\section}
{0pt}{0.05ex plus 0.01ex minus .05ex}{0.05ex plus .1ex}
\titlespacing*{\subsection}
{0pt}{0.05ex plus 0.01ex minus .05ex}{0.05ex plus .1ex}

\hypersetup{
	colorlinks = true, 
	linkcolor = blue,
	anchorcolor = red,
	citecolor = blue, 
	filecolor = red, 
}

\usepackage{algorithm,algorithmic}

\usepackage{graphicx}
\graphicspath{{./images/}}

\usepackage{titlesec}

\renewcommand\theparagraph{}

\titleformat*{\paragraph}{\bfseries}
\titleformat{\paragraph}[runin]
{\normalfont\normalsize\bfseries}
{\theparagraph}
{0em}
{}

\titlespacing*{\paragraph}
{0pt}
{3.25ex plus 1ex minus .2ex}
{1em}

\usepackage{pifont} % \ding{51} - i.e., tick mark

\usepackage{wasysym} % box with tick mark inside

\usepackage{amssymb}

\def\mytickbox2{$\text{\rlap{$\checkmark$}}\square$}
\def\mytickbox{$\text{\rlap{\ding{52}}}\square$}

\begin{document}
	%\begin{CJK*}{UTF8}{gbsn}
	
	\widowpenalties 1 1
	\raggedbottom
	\sloppy
	
	\title{\fontsize{16}{22}\selectfont \mytitle \\ \fontsize{13}{26}\selectfont \mysubtitle}

	%------------------------------------------------------------

	% \pdfminorversion=4
	
	\hypersetup{pdfinfo={
			Author		= {\myauthor},
			Title		= {\mytitle  \mysubtitle},
			Subject 		= {\mytitle  \mysubtitle},
			CreationDate = {D:20200220195600},
			Keywords 	= {\mykeywords},
	}}

	\author{\myauthor
		\IEEEcompsocitemizethanks{\IEEEcompsocthanksitem \myauthor \protect\\
			%\myemail 
		}% <-this % stops a space
		\thanks{}}
	
	%\author{Benjamin Kenwright,~\IEEEmembership{Member,~IEEE,}
		%       Idris Ibrahim
		%\IEEEcompsocitemizethanks{\IEEEcompsocthanksitem \myauthor \protect\\
			%\myemail 
			%}% <-this % stops a space
		%\thanks{}}

	\author{Benjamin Kenwright % <-this % stops a space
		\IEEEcompsocitemizethanks{\IEEEcompsocthanksitem % \protect\\
			Benjamin Kenwright ~\IEEEmembership{Senior,~IEEE,} bkenwright@ieee.org \\
			First draft: July 2020/Updated March 2023
		    }% <-this % stops a space
			
		\thanks{ }}

	%\author{ Benjamin Kenwright, Idris Ibrahim  }
	
	%\IEEEcompsocitemizethanks{
		%     \IEEEcompsocthanksitem Benjamin Kenwright
		%     \IEEEcompsocthanksitem Idris Ibrahim }
	% \thanks{ }

	\raggedbottom

	\markboth{ \mytitle (\myauthor) }%
	{ \mytitle }
	
	%------------------------------------------------------------

	%
	% Optional - teaser figure - top of the first page if you want
	% to show some statistics, figures, images, to capture the readers
	% interest
	
	%\teaser{
		%   \includegraphics[width=1.0\textwidth]{images/sampleteaser}
		%   \caption{Place a teaser image at the top of your report to show key examples of your work (e.g., multiple screenshots of the different test situations) - Every figure should have a caption and a description below it.  For example, each figure is labelled and explained: (a) shadows, (b) global illumination, (c) wireframe, and (e) motion blur.}

	%\bookmark[dest=tocpage,level=-1]{Contents}
	
	% \addcontentsline{toc}{section}{Chapter Virtual Reality: Where Have We Been? Where Are We Going?}

	%\pagenumbering{gobble}

	% \date{}
	
	% \maketitle
	
	\iffalse
	\section{abstract}
	\fi
	
	\vspace{-20pt}
	
	%% ------------------------------------------
	%% ------------------------------------------
	%%
	%% ABSTRACT ABSTRACT ABSTRACT ABSTRACT ABSTRACT
	%%
	%% ------------------------------------------
	%% ------------------------------------------
	
	%\section*{\centering\Large Abstract}
	
	\IEEEcompsoctitleabstractindextext{%
		\begin{abstract}
			\boldmath
			\myabstract
		\end{abstract}
		\begin{keywords}
			\mykeywords
	\end{keywords}}

	\maketitle
	
	\IEEEdisplaynotcompsoctitleabstractindextext
	\IEEEpeerreviewmaketitle
	
	%\setlength\cftparskip{-2pt}
	%\setlength\cftbeforesecskip{1pt}
	%\setlength\cftaftertoctitleskip{2pt}
	% \tableofcontents
	
	%-------------------------------------------------------------------------
	%-------------------------------------------------------------------------
	%-------------------------------------------------------------------------

	%\vspace{10pt}
	
	%\textbf{Keywords: \mykeywords}

%% ------------------------------------------
%% ------------------------------------------
%%
%% INTRODUCTION INTRODUCTION INTRODUCTION
%%
%% ------------------------------------------
%% ------------------------------------------

\section{Introduction}

Fourier transform (FT) is powerful mathematical tool used in many domains for analysing signal data \cite{kenwright2015quaternion}. 
In fact, the Fourier transform was probably one of the most important discoveries in the entire field of signal analysis.
However, while the Fourier transform was a game-changer that revolutionize science and engineering, its evolution and applications continue to evolve.
% It requires us to break out of our comfort zones and embrace a new way of thinking. Are you ready to take on the challenge and be a part of this groundbreaking change?
%
%
% Fourier analysis can also be applied to the characterization of spatial/frequency aspects.
%
%
%
Applying the Fourier transform to dual-quaternions unlocks a world of possibilities in robotics, computer graphics, and signal processing that were previously inaccessible (interpolate, analyze and extract new meanings from complex multicomponent data).
This approach has applications in signal processing, for example, filtering, analysis, reconstruction, and compression \cite{unuma1995fourier,lou2010example,shinya1998periodic}.
Dual-quaternions have been used to represent various multicomponent data types (rigid body transforms, motions, audio signals and image data).

The idea of computing the Fourier transform of a multicomponent signal has been realized previously.
However, there has been no definition of the Fourier transform applicable to dual-quaternion signals in a \textbf{holistic} manner.
In this paper, we present a method concerned with the computation of a unified single holistic Fourier transform of dual-quaternion multicomponent data. 
We apply this transform to motion filtering, such as, low pass smoothing or high pass filtering.
The application of quaternion Fourier transform (QFT) to signals is based on representing Euler joint changes using a quaternion discovered by Hamilton in 1843 \cite{shoemake1985animating}.  The first definition of a quaternion Fourier transform was that of Ell \cite{ell1993quaternion} and the first application of a quaternion Fourier transform was reported in 1996 for image processing using a discrete version of QFT \cite{guo2008spatio,sangwine1996fourier}.

While the \textbf{quaternion Fourier transform has gained much recognition, the concept has yet to be applied to dual-quaternion analysis and adaptation}.
We focus on demonstrating the enormous rewards of using the dual-quaternion Fourier transform (DQFT) for multicomponent data analysis (e.g., application in extracting and filtering signal components).

% \vspace{15pt}

The key contributions of this paper are:
(1) the application of dual-quaternion Fourier transform as a means of analysing multicomponent data signals;
(2) the holistic approach using Fourier transform of complex data as a frequency-domain representation based on the dual-quaternion numbers.

\section{Related Work}

The initial work on dual-quaternion FT (DQFT) was presented by Ding et al. \cite{ding2012approach}. Ding et al. presented multiple exciting insights around the application of DQFT for image processing techniques. % which is built upon a so called dual quaternion Fourier transform.
A lot of the concepts discussed in this article on using dual-quaternions are inspired from research around quaternion-based innovations.
For example, Kenwright \cite{kenwright2015quaternion}, work with quaternions and join signals was used to extend quaternion signals to dual-quaternion signals.  These dual-quaternion signals can then be transformed using the dual-quaternion Fourier Transform (as discussed in this article).
The idea of filtering signals to enhance the final pattern is not a new concept \cite{bruderlin1995motion,unuma1995fourier,lou2010example}, however, our work extends these techniques to incorporate the dual-quaternion Fourier transform and take advantage of the holistic nature approach.

While the Fourier transform has proven itself in many fields of engineering and computing for providing an un-cumbersome and efficient method of representing signal or functional information in the frequency domain, we hope to introduce the potential for expanding these concepts in the form of dual-quaternion Fourier transform for systems that need to amalgamate components into a unified form, such as, in motion analysis.

In the animation field, the Fourier series has been used on its own for classification and filtering.  One example, is Unuma et al. \cite{unuma1995fourier} who interpolated and transitioned between motion sequences using a Fourier series expansion. 
Similarly, quaternions are the de-facto method for efficiently representing unambiguous orientations that can be interpolated quickly and reliably \cite{shoemake1985animating}.
Wang et al. \cite{wang2006cartoon} used motion filtering techniques to add characteristics, such as, anticipation, follow-through, exaggeration and squash-and-stretch, which were not present in the
original input motion data. The filter adds a smoothed, inverted, or time shifted version of the `second derivative' of the signal back into the original signal.

Quaternion applications of Fourier transform for image processing \cite{bas2003color,sangwine1996fourier}, can be extended to dual-quaternions \cite{ding2012approach} or sound applications \cite{kenwright2023dual}. For a survey on dual-quaternions and their applications see Kenwright \cite{kenwright2023survey}.

Our solution combines the analogy of representing coupled data using dual-quaternion which can be transformed into the frequency domain using the dual-quaternion Fourier transform.

\section{Background}

\subsection{Quaternion Algebra} \label{sec:quaternions}

The concept of a quaternion was introduced by Sir. William Hamilton in 1843 \cite{shoemake1985animating}.
The quaternion is the generalization of a complex number.
A complex number has two components: a real and an imaginary part.
However, a quaternion has four components, i.e., one real and three imaginary parts.
A quaternion can be written in Cartesian form as:
\begin{equation}
	q = w + xi + yj + zk
\end{equation}

\noindent where $w$, $x$, $y$, and $z$ are real numbers and $i$, $j$, $k$ are complex operators which obey the following rules:
\begin{equation}  \label{eq:complex}
	\begin{alignedat}{3}
		& ij = k  \qquad jk && =i  \qquad  && ki =j\\
		& ji = -k \qquad kj && =-i \qquad  && ik=-j\\
		& i^2 = j^2=k^2     && =           && ijk = -1 
	\end{alignedat}
\end{equation}

The fundamental rules in Equation \ref{eq:complex}, emphasis that the multiplication of a quaternion is `not' commutative.  
Additional operations, include the quaternion conjugate ($q^*$):
\begin{equation}
	q^* = w - xi - yj - zk
\end{equation}

\noindent which is achieved by inverting the vector component.  The length (or modulus) of a quaternion is given by:
\begin{equation}
	|q| = \sqrt{w^2 + x^2 + y^2 + z^2}
\end{equation}

A quaternion with a zero real part is called a `pure quaternion', while a quaternion with a unit length (or modulus equal to one) is called a `unit quaternion':
\begin{equation}
	q = \frac{i + j + k}{\sqrt{3}} 
\end{equation}

The imaginary part of a quaternion has three components and may be associated with a 3-space vector.  For this reason, it is sometimes useful to consider a quaternion composition as a scalar part and a vector part, given by:
\begin{equation}  \label{eq:complex1}
	\begin{alignedat}{3}
		q &= w + xi + yj + zk\\
		&= s(q) + v(q)
	\end{alignedat}
\end{equation}

\noindent where $s(q)$ is the real part (i.e., $w$), and $v(q)$ is the vector part (i.e., $xi + yj + zk$).
\begin{equation}  \label{eq:complex2}
	\begin{alignedat}{3}
		s(q) &= w\\
		v(q) &= xi + yj + zk
	\end{alignedat}
\end{equation}

A 3-dimensional joint with three degrees of freedom has three Euler angles (i.e., three components, $x$, $y$, and $z$) that we represent using a `pure quaternion'.
For the joint signal in xyz local space, the three imaginary parts of a pure quaternion can be used to represent the $x$, $y$, and $z$ axis of rotation component.
A joint from an animation at time $t$ can be represented as:
\begin{equation}
	f(t) = x(t)i + y(t)j + z(t)k
\end{equation}

\noindent where $x(t)$, $y(t)$, and $z(t)$ are the $x$, $y$, and $z$ Euler angle displacements at time $t$.
Using a quaternion to represent the xyz joint space transforms, the three channel components are processed equally in operations, such as, multiplication.
The advantage of using a quaternion based operation to manipulate the angular information for a joint is that we do not process each axis independently, but rather, treat each joint as a unified unit.
We believe, by using a quaternion as a unified unit, improves motion information accuracy because a joint is treated as an entity.

% ref: https://twitter.com/SamuelGWalters/status/1425726883283406850
\subsection{Quaternion Exponentials and Logarithms}
A quaternion (q) has an exponential defined as:

\begin{equation} \label{eq:quatlog}
	e^q = cos(|q|) + sin(|q|) \frac{q}{|q|}
\end{equation}

All non-zero quaternions have logarithms and polar form values. For $q$ in Equation \ref{eq:quatlog}, $q^2=-|q|^2$, which says that $q/|q|$ behaves like the complex $i$. If $q=0$, the right side of the equation is understood to be $1$.

\subsection*{Dual-Quaternion Algebra}
For completeness, we include the mathematical details for essential dual-quaternion operations.

\paragraph{Dual-Quaternion Arithmetic Operations}

The elementary arithmetic operations necessary for us to use dual-quaternions are given below:

\begin{itemize}
	
	\item \textbf{dual-quaternion}: $\zeta = q_r + q_d \varepsilon$
	
	\item \textbf{scalar multiplication}:  $s \zeta = s q_r + s q_d \varepsilon$
	
	\item \textbf{addition}: $\zeta_1 + \zeta_2 = q_{r1} + q_{r2} + (q_{d1}+q_{d2}) \varepsilon$
	
	\item \textbf{multiplication}: $\zeta_1 \zeta_2 = q_{r1} q_{r2} + (q_{r1}q_{d2} + q_{d1}q_{r2}) \varepsilon$
	
	\item \textbf{conjugate}:  $\zeta^* = q_r^* + q_d^* \varepsilon$
	
	\item \textbf{magnitude}: $||\zeta|| = \zeta \zeta^{*}$
	
\end{itemize}
\noindent where $q_r$ and $q_d$ indicate the real and dual part of a dual-quaternion (we use a `$q$' to emphasis that they are quaternions).

\subsection{Dual-Quaternion Exponentials and Logarithms}

The dual-quaternion exponential function is a mathematical function that takes a dual-quaternion as an input and returns another dual quaternion as the output. 
The form for dual-quaternion exponential is well established \cite{selig2010exponential,kenwright2022inverse}, additional variations have been presented \cite{dantam2021robust} to addresses the zero-angle singularity in the dual-quaternion exponential and logarithm (practical necessity if dual-quaternions are to be for Fourier transformations).
% l

%A dual-quaternion function of time $t$ can be represented as:
%
%\begin{equation}
%	f(t) = x(t)i + y(t)j + z(t)k
%\end{equation}

For a beginners introduction to dual-quaternions and a comparison of alternative methods (e.g., matrices and Euler angles) and how to go about implementing a straightforward library we refer the reader to the paper by Kenwright \cite{kenwright2012beginners}.

\paragraph{Dual-Quaternion Vector Transformation}
A dual-quaternion is able to transform a 3D vector coordinate as shown in Equation \ref{eq:dqtransf}. Note that for a unit dual-quaternion the inverse is the same as the conjugate \footnote{
	Just to note, there are three definitions for the conjugate of a dual quaternion: 
	1: $\zeta^* = q_r^* + q_d^* \varepsilon$,
	2: $\zeta^* = q_r - q_d  \varepsilon$, and
	3: $\zeta^* = q_r^* - q_d^* \varepsilon$.
	For transforming a point using Equation \ref{eq:dqtransf}, you would use the 2nd conjugate variation.
	% ref: https://www.euclideanspace.com/maths/algebra/realNormedAlgebra/other/dualQuaternion/functions/index.htm
}.

\begin{equation}
	p' = \hat{\zeta} p \hat{\zeta}^{-1}
	\label{eq:dqtransf}
\end{equation}

\noindent where $\hat{\zeta}$ is a unit dual-quaternion representing the transform, $\hat{\zeta}^{-1}$ is the inverse of the unit dual-quaternion transform.
$p$ and $p'$ are the dual-quaternions holding 3D vector coordinate to before and after the transformation (i.e., $p=(1,0,0,0) + \epsilon (0,v_x,v_y,v_z)$ )).

\paragraph{{Pl{\"u}cker Coordinates}}
Pl{\"u}cker coordinates are used to create Screw coordinates which are an essential technique of representing lines. We need the Screw coordinates so that we can re-write dual-quaternions in a more elegant form to aid us in \textbf{formulating a neater and less complex interpolation method} that is comparable with spherical linear interpolation for classical quaternions.

The Definition of Plu\"cker Coordinates:
\begin{itemize}
	\item $\vec{p}$ is a point anywhere on a given line
	\item $\vec{l}$ is the direction vector
	\item $\vec{m} = \vec{p} \times \vec{l}$ is the moment vector
	\item $(\vec{l},\vec{m})$ are the six Plu\"cker coordinate
\end{itemize}

We can convert the eight dual-quaternion parameters to an equivalent set of eight screw coordinates and vice-versa.
The definition of the parameters are given below in Equation \ref{eq:dqscrew}:

\begin{equation}
	\begin{alignedat}{4}
		\text{screw parameters} &= (\theta,d,\vec{l},\vec{m}) \\
		\text{dual-quaternion}  &= q_r + \epsilon q_d \\
		&= (w_r+\vec{v}_r) + \epsilon(w_d + \vec{v}_d) 
	\end{alignedat}
	\label{eq:dqscrew}
\end{equation}

\noindent where in addition to $\vec{l}$ representing the vector line direction and $\vec{m}$ the line moment, we also have $d$
representing the translation along the axis (i.e., pitch) and the angle of rotation $\theta$.

\subsection{Discrete Dual-Quaternion Fourier Transform (DQFT)}
The mathematical principles introduced in Section \ref{sec:quaternions} define the fundamental concepts, such as, multiplication and exponential, that are used in the Discrete Quaternion Fourier Transform (DQFT) (see Sangwine \cite{ell2014quaternion,said2008fast}).
A crucial factor when using DQFT is dual-quaternion multiplication is non-commutative, hence, there are two different types of DQFT operation. %
The two types are defined as, the right side DQFT and the left side DQFT.  The two DQFT types are given below in Equations \ref{eq:dqftR} and \ref{eq:dqftL}:

\noindent Right sided DQFT:
\begin{equation}  \label{eq:dqftR}
	\begin{alignedat}{3}
		F_{R}(t) = \frac{1}{\sqrt{M}} \sum_{x=0}^{M-1}  f(x) \; e^{-\mu 2 \pi \left( \frac{x t}{M} \right)}
	\end{alignedat}
\end{equation}

\noindent Left sided DQFT:
\begin{equation}  \label{eq:dqftL}
	\begin{alignedat}{3}
		F_{L}(t) = \frac{1}{\sqrt{M}} \sum_{x=0}^{M-1}  e^{-\mu 2 \pi \left( \frac{x t}{M} \right)} \; f(x)
	\end{alignedat}
\end{equation}

%\noindent Right-Left sided DQFT:
%\begin{equation}  \label{eq:dqftRL}
%\begin{alignedat}{3}
%F_{L}(t_0,t_1) = \frac{1}{\sqrt{MN}} 
%\sum_{x=0}^{M-1} \sum_{y=0}^{N-1} 
%e^{-\mu 2 \pi \left( \frac{x t_0}{M} \right)} 
%\; f(x,y) \;
%e^{-\mu 2 \pi \left( \frac{y t_1}{N} \right)} 
%\end{alignedat}
%\end{equation}

Likewise, we are able to define the `Inverse' Dual-Quaternion Fourier Transform (IDQFT) for the two types of DQFT respectively, as given below in Equations \ref{eq:idqftR} and \ref{eq:idqftL}:

\noindent Right sided IDQFT:
\begin{equation}  \label{eq:idqftR}
	\begin{alignedat}{3}
		f_{R}(t) = \frac{1}{\sqrt{M}} \sum_{x=0}^{M-1}  f(x) \; e^{\mu 2 \pi \left( \frac{x t}{M} \right)}
	\end{alignedat}
\end{equation}

\noindent Left sided IDQFT:
\begin{equation}  \label{eq:idqftL}
	\begin{alignedat}{3}
		f_{L}(t) = \frac{1}{\sqrt{M}} \sum_{x=0}^{M-1} e^{\mu 2 \pi \left( \frac{x t}{M} \right)} \; f(x)
	\end{alignedat}
\end{equation}

%\noindent Right-Left sided IDQFT:
%\begin{equation}  \label{eq:idqftRL}
%\begin{alignedat}{3}
%f_{L}(t_0,t_1) = \frac{1}{\sqrt{MN}} 
%\sum_{x=0}^{M-1} \sum_{y=0}^{N-1} 
%e^{\mu 2 \pi \left( \frac{x t_0}{M} \right)} 
%\; f(x,y) \;
%e^{\mu 2 \pi \left( \frac{y t_1}{N} \right)} 
%\end{alignedat}
%\end{equation}

% \noindent where $\mu$ is any pure unit dual-quaternion and determines a direction in three-dimensional space.  This maybe seen as a control factor - typically, a suitable value is chosen to connect all the angles (i.e., for coupled influence).

%\figuremacroW
%{freqrange}
%{Animation Frequencies}
%{Joint motions are essentially signals that occupy specific frequency ranges.}
%{1.0}

\section{Applications of Dual-Quaternion Fourier Transform}

The Dual-Quaternion Fourier transform has a wide range of applications, including:

\begin{enumerate}
\item Signal Processing: Dual-Quaternion Fourier transform is used to analyze and process signals in areas such as audio and video processing, image processing, and telecommunications. It is used to convert a time-domain signal to a frequency-domain signal, which makes it easier to analyze and process.

\item Data Compression: Dual-Quaternion Fourier transform is used in data compression techniques such as JPEG and MPEG, where it is used to reduce the amount of data required to represent an image or video by eliminating the high-frequency components that are not perceptible to the human eye.

\item Optics: It is used in the study of optics, including the analysis of diffraction patterns and the design of optical systems such as lenses and telescopes.

\item Quantum Mechanics: Dual-Quaternion Fourier transform is used in quantum mechanics to describe the wave functions of particles, which allows us to calculate the probability of finding a particle in a particular location.

\item Medical Imaging: Dual-Quaternion Fourier transform is used in medical imaging techniques such as MRI and CT scans to reconstruct images of the human body from signals generated by the scanning process.

\item Speech Recognition: Dual-Quaternion Fourier transform is used in speech recognition systems to convert the speech signal from the time-domain to the frequency-domain, which makes it easier to identify speech patterns and recognize words.

\item Climate Analysis: Dual-Quaternion Fourier transform is used in climate analysis to analyze temperature and other climate data over time, which helps in understanding climate patterns and making predictions about future climate trends.

\end{enumerate}

These are just a few examples of the wide range of applications of the Dual-Quaternion Fourier transform in various fields. However, as there is a large body of research around spatial transformation signals (e.g., animation and motion of rotations and translations) \cite{kenwright2023survey}, we present additional information on `filtering' using the DQFT.

\subsection{DQFT Filtering}

The dual-quaternion Fourier transform (DQFT) takes the signal variables at time $(t)$ performs a filtering operation in the frequency domain.
Importantly, the dual-quaternion convolution operation in the spatial domain corresponds to the product operation in the frequency domain \cite{pei2001efficient}.
After filtering the data in the frequency domain, we are able to recover the modified signal by taking the inverse dual-quaternion Fourier transform.
%

%\figuremacroW
%{filtertypes}
%{Filter}
%{Gradual and short filter types.}
%{0.9}

\paragraph{DQFT Low-Pass \& DQFT High-Pass}
The low-pass filter attenuates a specified range of high frequency components in the dual-quaternion Fourier transform signal.
The DQFT low-pass filter provide an extremely powerful mathematical tool for separating a signal in its different components.  Once we are able to decompose and filter signal components it makes the problem much more simple to analyze.
A high-pass filter attenuates low-frequency components without disturbing high frequency information in the dual-quaternion Fourier transform.
The high-pass filter is essentially the reverse of the low-pass filter, which we invert to accomplish the desired operation.

\paragraph{DQFT Range-Pass}
Combining both high and low-pass filters, we are able to filter `targeted' signals in the dual-quaternion Fourier transform for either attenuation or extraction.  This might include, rhythmic noise that we want to collect for re-injection into other sources for data classification by identifying key frequencies that make up specific patterns.
Similarly, by comparing signals and identifying similarities using the dual-quaternion Fourier transform, we are able to extract and differentiate common and unique signal characteristics.
This work feeds into a number of key fields, such as, medical analysis, sound, procedural animations and robotics.

%\figuremacroW
%{freqdomain}
%{Spatial-Frequency Domain}
%{Taking joint signals from the spatial domain to the frequency domain we are able to decompose complex actions into simpler pieces.}
%{1.0}

%\section{Experimental Results}
%We imported and applied DQFT filtering operations to a variety of motion capture actions (e.g., walking, falling, and dancing).
%For joints with fewer than three degrees of freedom (e.g., knee joint), the quaternion was padded with a constant %value to formulate three Euler angles for the DQFT operation.  
%Since the padded Euler values were constant they did not influence or affect the joint signal that was changing.
%Of course, due to the articulated character hierarchical topology, signal modifications would influence children limbs and ultimately end-effector positions and orientations, for instance, feet and hands.
%We only focused on the filtering but the end-effectors would need re-targeting either using inverse kinematic methods or artistic intervention (see Figure \ref{fig:anklesignals} and \ref{fig:freqdomain}).
%Figure \ref{fig:filteringresults} and Figure \ref{fig:filteringresults3} show examples of motion filter simulations. 

\paragraph{Benefits of Dual-Quaternion \& Fourier Transforms}

\begin{itemize}
	\item Helps reduce sensitivity to noise, generally, the Fourier Transform is sensitive to noise, meaning that even small amounts of noise can distort the spectral content of the signal. This sensitivity can make it difficult to accurately analyze signals with low signal-to-noise ratios. Due to the coupled multicomponent data configuration of dual-quaternions this sensitive to noise is reduced.
	
	\item Enhance time-frequency resolution - by default the Fourier Transform assumes that the signal is stationary, meaning that it does not change over time. This assumption limits the ability of the Fourier Transform to capture non-stationary signals, such as signals with time-varying frequency components or signals with transient features. % The resolution of the Fourier Transform in the time and frequency domains is also limited, meaning that it cannot precisely determine the timing of events or the frequency content of signals with closely spaced spectral components.
	Coupled multicomponent data of the dual-quaternion helps improve the resolution for determining timing of events or the frequency content of signals with closely spaced spectral components.
	
	% Assumes Infinitely Long Signals: The Fourier Transform assumes that the signal extends infinitely in both the positive and negative time directions. This assumption is not valid for finite-length signals, and can result in spectral leakage and other artifacts.
	
	\item Improves spatial resolution, for example, in image analysis and other applications, the Fourier Transform is often used to analyze the frequency content of images. However, the Fourier Transform by default limited spatial resolution, meaning that it cannot precisely localize frequency components in space. This limitation can make it difficult to analyze complex images with multiple features or regions of interest.

	\item Dual-quaternion representation of the data can help link related information to improve interpretability.  For instance, the basic Fourier Transform provides a mathematical representation of the frequency content of a signal, but it does not provide direct interpretations of the signal. Interpretation of the Fourier Transform requires additional knowledge of the signal and its properties, which can be difficult to obtain in some cases using other methods).

\end{itemize}

\section{Discussion}
Dual-quaternions have become a popular mathematical tool for multicomponent data, offering a efficient, compact and unified solution for manipulating, concatenating and interpolating data. % 
%However, the dual-quaternion Fourier transform in this paper does not use the joint quaternion orientation transforms directly, in retrospect, the Euler angles are transformed to the QFT space for filtering.
However, as we have presented in this paper, the dual-quaternion representation can be combined with the Fourier transform to reap additional opportunities.
%
%We suggest a \textbf{holistic Fourier transform} of the joints to yield a single frequency-domain representation based on the quaternion Fourier coefficients.
For example, combining the dual-quaternions with the Fourier transform addresses a number of limitations (including sensitivity to noise and improved interpretability). % . this addresses a number of limitat opens the door to new types of motion filtering techniques.
%
% Any adaptation of the original motion, such as, removal or modification of the joint signals, will cause the overall animation to change.  Due to the `hierarchical' nature of the model (i.e., skeleton of connected limbs), that propagate transforms.  Modifying a joint's local coordinates in the hierarchy would cause inherited limbs world coordinates to change.  To ensure `end-effectors' (e.g., hands, and feet) reach their in-tented target the final motion is blended with an inverse kinematic model.
%
% Since the underlying transforms are `Euler' angles, which suffers from the famous `gimbals lock'.  The xyz Euler angles are used compared to other methods to keep the underlying linear space.
%
% A variety of motion capture file formats store the transforms as Euler angles for each axis (i.e., xyz) allowing the data to be loaded and filtered easily without pre-processing.
%
The filters are linear at their core, however, dual-quaternions are able to embrace highly non-linear forms  constructing filters that are more inapt. 

\section{Conclusion}
This article has shown that dual-quaternion Fourier transform (DQFT) is well suited for multiple areas including adapting and filtering the spectral content of signal data.
For example, six-dimensional joint motions represented as dual-quaternions can be transferred efficiently to the frequency domain (represented as a dual-quaternion frequency signal) for processing techniques, such as, filtering.
Filtering in the dual-quaternion frequency domain has the advantage that the motion data is processed as a whole unit rather than dealing with the individual channels separately.
We believe more accurate signal information can be preserved this way, since the signal channels are processed as a single unit.
This article has also presented other applications and research beyond filtering.
All things considered, we like to point out that, the dual-quaternion Fourier transform is not limited to transform data (e.g., rotations and translations), but can also be applied to other forms and processes, such as, pattern recognition and data compression.
Our work serves to complement the many interactive tools available for editing and constructing signal using dual-quaternions as the base-format.

\bibliographystyle{plain}

\let\oldthebibliography=\thebibliography
\let\endoldthebibliography=\endthebibliography

\renewenvironment{thebibliography}[1]{%
    \begin{oldthebibliography}{#1}%
      \setlength{\parskip}{0ex}%
      \setlength{\itemsep}{0.5ex}%
  }%
  {%
    \end{oldthebibliography}%
  }

\fontsize{8}{8.0}\selectfont

\bibliography{paper}

%\end{CJK*}
\end{document}